\documentclass[12pt]{amsart}
\usepackage{geometry}                
\usepackage{graphicx}
\usepackage{amssymb}
\usepackage{epstopdf}
\usepackage{enumerate}

\theoremstyle{plain} 
\newtheorem{theorem}{Theorem} 
\newtheorem{proposition}{Proposition} 
\newtheorem{lemma}{Lemma} 
\newtheorem*{corollary}{Corollary}

\theoremstyle{definition} 
\newtheorem{definition}{Definition} 
\newtheorem{conjecture}{Conjecture} 
\newtheorem{claim}{Claim}
\newtheorem*{example}{Example}

\theoremstyle{remark}

\title{Expanding Polynomials over the rationals}

\author[1]{J\'ozsef Solymosi}

\address{Department of Mathematics, University of British Columbia \\ 1984 Mathematics Road, Vancouver, BC, V6T 1Z4, Canada}
\email{solymosi@math.ubc.ca}
\thanks{Research was supported by NSERC and OTKA grants}
\begin{document}

\begin{abstract} Let $F(x,y)$ be a polynomial over the rationals. We show that if $F$ is not an expander (over the rationals) then it has a special multiplicative or additive form. For example if $F$ is a homogeneous non-expander polynomial then $F(x,y)=c(x+ay)^\alpha$ or $F(x,y)=c(xy)^\alpha .$ This is an extension of an earlier result of Elekes and R\'onyai who described the structure of two-variate polynomials which are not expanders over the reals.
\end{abstract}

\maketitle

\section{Introduction}
Let $f: {\mathbb{C}\times \mathbb{C}\rightarrow \mathbb{C}}$ be a two-variable polynomial over the complex numbers. For two sets $A,B\subset {\mathbb{C}}$ we denote the range of $f$ over $A\times B$ by $f(A,B).$ 

\[
f(A,B)=\{ f(a,b)\; |\; a\in A, b\in B\}.
\]

In this paper we are interested about polynomials for which $|f(A,B)|\ll \max{(|A|,|B|)}$ for some large sets of rational numbers $A,B \subset {\mathbb{Q}}.$ To formulate our statements we use the following definition.

\begin{definition}Let $f$ be a $k$-variable polynomial over some base ring $R.$ We say that $f$ is a  {\em non-expander} over a domain $D\subset R$ if there is constant $c>0$ so that for any integer $n$ there are sets $A_1,A_2,\ldots ,A_k\subset D$ with $|A_1|=|A_2|=\ldots=|A_k|\geq n$ and $|f(A_1,A_2,\ldots ,A_k)|\leq c|A|.$ Otherwise polynomial $f$ is an {\em expander}. 
\end{definition}

Elekes and R\'onyai  described the structure of two-variate polynomials which are non-expanders over the reals \cite{ER}. They proved that a polynomial, $F(x,y),$ is non-expander over the reals if and only if
\begin{equation}\label{ER}
F(x,y)=f(g(x)+h(y))\;\; or \;\; F(x,y)=f(g(x)h(y))
\end{equation}
for some (single variable) polynomials $f,g,$ and $h.$ 

Our plan -- which we can not complete here without additional assumptions -- is to show that  a polynomial, $F(x,y),$ is non-expander over the rationals if and only if \begin{equation}\label{main}
F(x,y)=f(x+ay)\;\; or \;\;F(x,y)=f((x+a)^{\alpha}(y+b)^{\beta})
\end{equation} 
for some (single variable) polynomial $f,$ positive integers $\alpha,\beta,$ and rational numbers $a,b.$ 

\subsection{Previous results}
The result of Elekes and R\'onyai was extended by Elekes and Szab\'o \cite{ESz} for implicit polynomials providing a larger scale of applications. The typical applications of the mentioned results are arising from Erd\H os-type problems in discrete geometry. It was Elekes who pioneered the use of such algebraic methods in discrete geometry. (See Elekes' survey \cite{EL} on the subject.) 

In a recent paper Tao \cite{Ta} proved an Elekes-R\'onyai type theorem for the finite field case. There were many interesting and important predecessors of his result. We refer the reader to \cite{Bo} and \cite{Ta} for more details about the finite field case. We only mention here one example since it is about a three-variable polynomial which we will consider later. Shkredov proved  in \cite{Sh}  that $f(x_1,x_2,x_3)=x_1^2+x_1x_2+x_3$ is an expander over finite fields. Using a generalization of the Elekes-R\'onyai theorem for three variables, proved in \cite{SSZ}, we chacterize  three-variable polynomials which are non-expanders over the rationals.     

In this paper we are interested about polynomials which are expanders over the rationals. In this direction Chang proved in \cite{Ch} that $f(x,y)=x^2+y^2$ is an expander over the rationals. Note that this polynomial is a non-expander over the reals. Selecting $A=\{1,\sqrt{2},\sqrt{3},\ldots ,\sqrt{n}\}$ one can see that $|f(A,A)|\leq 2|A|.$ 

\section{Results}
Polynomials which are non-ex\-pan\-ders over the rationals are obviously non-ex\-pan\-ders over the reals.  By the result of Elekes and R\'onyai it is known that such polynomials  have a multiplicative or additive structure as in (\ref{ER}). In order to get a more restrictive form, as in (\ref{main}), one should analyze the inner polynomials $h(x)$ and $g(y).$ If $F(x,y)=f(g(x)+h(y))$ then $|g(A)+h(B)|\leq d|F(A,B)|\leq d|g(A)+h(B)|$ where $d$ is the degree of $F.$ Thus a Pl\"unecke-Ruzsa type inequalty, \cite{PR}, \cite{Pl}, implies that if $|F(A,B)|\leq Cn$ for some large $|A|=|B|=n$ then $|g(A)+g(A)|\leq C'n$ where $C'$ depends on $C$ only. Similarly, if $F(x,y)=f(g(x)h(y))$ then $|g(A)g(A)|\leq C'n$. One can apply Freiman's theorem, \cite{Fr}, \cite{Ru},  to see that $g(A)$ is a dense subset of a $D=D(C)$-dimensional generalized arithmetic or geometric progression,
\[ GAP_D=\sum_{i=1}^D AP_{n_i}\;\;\; or \;\;\; GGP_D=\prod_{i=1}^D GP_{n_i},\]
where $AP_{n_i}$ is a set of the elements of an $n_i$-term arithmetic progression and 0, and $GP_{n_i}$ is a set of the elements of an $n_i$ element geometric progression and 1. For more details about generalized arithmetic progressions and Freiman's theorem we refer to \cite{TV}. Freiman's theorem states that not only the dimension of the GAP or GGP is bounded (this is often called Freiman's Lemma and it is relatively easy to prove) but $g(A)$ is dense in the generalized progression, i.e. $\prod_{i=1}^Dn_i\leq C''n$ for some $C''=C''(C).$   We will use the stronger version for the additive case in \ref{main}. For the multiplicative case the Freiman Lemma is sufficient. 

\subsection{The multiplicative case} 

Here we suppose that  $F(x,y)=f(g(x)h(y))$ and $|g(A)g(A)|\leq C'n$. The set  $g(A)$ is a (dense) subset of a $D=D(C)$-dimensional generalized geometric progression,
$GGP_D=\prod_{i=1}^D GP_{n_i}.$ In particular, the range of $g(x)$ over the rationals has an $n$-element intersection with a multiplicative group of rank at most $D.$ As we will see, this is only possible if $g(x)=(x+a)^\alpha$ for some rational $a$ and integer $\alpha.$

\begin{lemma}\label{multi}
For any polynomial over the rationals, $g(x),$ there is a bound, $b=b(g,r),$ such that for any complex multiplicative group $G$ of rank $r$ the size of the intersection, $|g({\mathbb Q})\cap G|,$ is at most $b$ unless
$g(x)=(x+a)^\alpha$ for some rational $a$ and integer $\alpha.$
\end{lemma}

From Lemma \ref{multi} and the previous observations it follows that if $F(x,y)$ is a non-expander over the rationals and it has the multiplicative form in (\ref{ER}) then it has the more restricted form in (\ref{main}). 

\medskip
\noindent
{\em Proof of Lemma \ref{multi}:} Write 
$$g(x)=\prod_{i=1}^v(x-x_i)^{\alpha_i}$$ 
where $x_i$-s are the roots of $g$. If a natural number $w> deg(g)=\sum_{i=1}^v\alpha_i$ is relative prime to $w-deg(g)$ and to all $\alpha_i$-s then the algebraic curve 
\begin{equation}\label{curve}
{\mathcal{C}}:\;\;g(x)=y^w
\end{equation}
has genus 
$$\frac{(v-1)(w-1)}{2}.$$ If $v>1$ and $w$ is prime so that $w\geq deg(g)\geq 5$ then $\mathcal{C}$ has genus at least two. By Faltings' theorem \cite{Fa} it has finitely many rational points. 

Now let us consider the intersection $g({\mathbb Q})\cap G.$ $G$ is a finitely generated Abelian group and its rational elements form a subgroup. This subgroup is also finitely generated and its rank is not larger than $rank(G).$ Therefore in Lemma \ref{multi} we can suppose without loss of generality that $G$ is a multiplicative subgroup of ${\mathbb Q}.$   Suppose that $G$ is generated by elements $a_1, a_2,\ldots ,a_s.$  There is a natural surjection of $\mathbb{Z}^s$ onto $G$
\[ \xi(\beta_1,\ldots ,\beta_s) \mapsto a_1^{\beta_1}a_2^{\beta_2}\cdots a_s^{\beta_s}.\]
We are interested about the the set $$S=\xi^{-1}(G\cap g({\mathbb{Q}})).$$
At least $|S|/w^s$ of the elements of $S$ are in the same congruence class modulo ${w}.$ The representative of this popular congruence class is denoted by $(w_1,w_2,\ldots ,w_s)$ where $0\leq w_i < w.$ The corresponding elements of $S$ raise rational points on the curve
\begin{equation}
\mathcal{C'}: \;\; g(x)=a_1^{w_1}a_2^{w_2}\cdots a_s^{w_s}y^w.
\end{equation} 
This is a twist of curve $\mathcal{C}$ in (\ref{curve}). According to a result of Silverman \cite{Si}, which is an extension of Faltings'  theorem, the number of points on $\mathcal{C'}$ is uniformly bounded by a constant depending on $\mathcal{C}$ only.  This also gives a uniform bound on $|S|,$ so $v=1$ which proves the Lemma. 

\medskip
\qed

\medskip
\noindent
The bound in Lemma \ref{multi} depends on $g(x)$ but there might be a uniform bound depending on the degree of $g(x)$ only.  Such bound would follow from the Uniformity Conjecture, see in \cite{CHM} and \cite{AV}. The conjecture states that there is a uniform bound on the rational points of a curve over the rationals depending on the genus of the curve only. 

\begin{conjecture}\label{unif-multi}
There is a function $b(d,r)$ so that the following holds. For any degree $d$ polynomial over the rationals, $g(x),$ and for any complex multiplicative group, $G,$ of rank $r$ the size of the intersection $g({\mathbb Q})\cap G$ is uniformly bounded by $b(d,r)$ unless $g(x)=(x+a)^\alpha$ for some rational $a$ and integer $\alpha.$\end{conjecture}

\subsection{The additive case}
We are going to use the claim that sets with small sumset contain a long arithmetic progression. 
\begin{proposition}\label{AP}
 For every integer $\ell$ and constant $D$ there is a threshold $n_0=n_0(\ell,D)$ such that if $|g(A)+h(B)+u(C)|\leq Dn$ where $|A|=|B|=|C|=n\geq n_0$ then $g(A)$ contains an $\ell$-term arithmetic progression.
\end{proposition}
This proposition follows from a combination of Freiman's theorem and Szemer\'edi's theorem about long arithmetic progressions, but one can see it directly using the Hypergraph Removal Lemma as it was shown in \cite{So}.

Because of Proposition \ref{AP}, in order to show that non-expander polynomials have the restricted multiplicative or additive form in (\ref{main}), it would be enough to prove the following statement which we state as a conjecture. 
\begin{conjecture}\label{APconj}
 For any rational polynomial $g(x)$ there is a bound $b$ so the range over the rationals, $g(\mathbb{Q}),$ contains no arithmetic progression of length $b.$ 
 \end{conjecture}
 Unfortunately we were unable to prove the conjecture above. It would follow from the Uniformity Conjecture mentioned earlier.
 \begin{claim}
Conjecture \ref{unif-multi} implies  Conjecture \ref{APconj}.
\end{claim}
\noindent
{\em Proof:} Let us suppose that $g(\mathbb{Q})$ has a large (compare to $deg(g)$), say $N$-term arithmetic progression $a+id.$ If $x_0$ denotes the rational number for which $g(x_0)=a$ then the polynomial $h_{a}(x)=g(x+x_0)-g(x_0)$ satisfies  the $h_{a}(x_i)=id$ equations for some $x_i$ $(1\leq i\leq N)$ rational numbers. Note that it is not possible that $h_{a}(x)=(x-\rho)^\alpha$ for some rational number $\rho$ if $\alpha \geq 2$ and $N\geq 4$ since the set of squares contains no four-term arithmetic progression and any set of larger powers $\{j^\alpha\}_{j=1}^\infty$ contains no three-term arithmetic progression. Even if $\rho$ is not rational, the extension field of $\mathbb{Q}$ with $\rho$ has similar properties, there are no long arithmetic progressions among the perfect squares over this field \cite{Xa} and there are no long arithmetic progressions for larger $\alpha$-s, see e.g. in \cite{Si}. Therefore we can suppose that $h_{a}(x)$ has at least two distinct roots. It also contains a long geometric progression $h_{a}(x_{2^j})=d2^j$ for $2^j\leq N$ which for large enough $N$ contradicts with Conjecture \ref{unif-multi}. 

\qed

\medskip

There are various special cases of the Uniformity Conjecture which would imply the additive part of formula (\ref{main}). We were unable to prove even the integer case, that for a given integer polynomial, $g(x),$ which has  degree at least two,  the range, $g(\mathbb{Z}),$ contains no arbitrary long arithmetic progressions.  The following is maybe the weakest form of a uniformity conjecture which would imply it.
\begin{conjecture}
For any polynomial over the integers $g$ there is a bound $b$ so that for any $a\in \mathbb{Z}$ the range of the polynomial $h_a(x)=g(x+a)-g(a)$ contains no $b$-term geometric progression. ($GP_b\notin h_a(\mathbb{Z}).)$ 
\end{conjecture}
There are extensions of the Elekes-R\'onyai theorem for three-variable polynomials and for sets $A,B,C$ having different sizes. In \cite{SSZ} we proved that if $F(x,y,z)$ is a non-expander over the reals then 
\begin{equation}\label{ER3}
F(x,y,z)=f(g(x)+h(y)+u(z))\;\; or \;\; F(x,y,z)=f(g(x)h(y)u(z))
\end{equation}
for some polynomials $f,g,h,$ and $u.$
Under the hypothesis that the Uniformity Conjecture holds we can conclude that if $F(x,y,z)$ is not expander over the rationals then 
\begin{equation}\label{3var}
F(x,y,z)=f(x+ay+bz)\;\; or \;\;F(x,y)=f((x+a)^{\alpha}(y+b)^{\beta}(z+c)^\gamma)
\end{equation} 
for some (single variable) polynomial $f,$ positive integers $\alpha,\beta,\gamma$ and rational numbers $a,b,c.$ 

\section{The unconditional case}
By putting together the previous arguments we have the following theorem, which does not require the Uniformity Conjecture.
\begin{theorem}\label{vegso}
If $F(x,y,z)$ is not expander over the rationals then 
\begin{enumerate}[i]
\item , $F(x,y,z)=f(x+ay+bz)$ or 
\item , $F(x,y,z)=f(g(x)+h(y)+u(z))$ or 
\item , $F(x,y,z)=f((x+a)^{\alpha}(y+b)^{\beta}(z+c)^\gamma)$
\end{enumerate}
for some (single variable) polynomial $f,$ polynomials $g, h, u$ without constant terms which all have degree at least three and have at least two distinct roots, positive integers $\alpha,\beta,\gamma$ and rational numbers $a,b,c.$ 
\end{theorem}
We believe that case (ii) is not possible, however even this theorem is useful to show if a polynomial is an expander over the rationals or not. The situation is much simpler with homogeneous polynomials.

\begin{corollary}
Every homogeneous polynomial $F(x,y,z)$ is an expander over the rationals unless it has the form  
\begin{enumerate}[i]
\item , $F(x,y,z)=a(x+by+cz)^\alpha$ or 
\item , $F(x,y,z)=a(xyz)^\alpha$
\end{enumerate}
for some rationals $a,b,c$ and integer $\alpha.$
\end{corollary}
We leave the simple proof -- which is a case analysis based on Theorem \ref{vegso} -- to the reader.

\subsection{Examples} As we mentioned earlier, Shkredov showed  in \cite{Sh}  that $$F(x,y,z)=x^2+xy+z$$ is an expander over finite fields (note that the definition of expanders is somewhat different over finite fields.) We show that $F$ is also expander over the rationals. By Theorem \ref{vegso} we know that the forms in (ii) or (iii) are not possible since the degree of $F$ is two, and in both cases (ii) and (iii) the polynomial has degree at least three. Case (i) is not possible either because $F$ is not a quadratic polynomial of $x+ay+cz.$  

\medskip

Let us see some more complex examples. The main purpose of these examples is to illustrate how to decide if a polynomial is an expander over the rationals or not by using Theorem \ref{vegso}. 
\begin{example}
The polynomial $$F(x,y,z)=x^2y^2z+x^2y^2+x^2z+x^2-y^2z-y^2-z-1$$
is an expander over the rationals (and a non-expander over the reals).
\end{example}

\medskip
\noindent
{\em Proof:} This is a degree five polynomial, so if it had the form (i) in Theorem \ref{vegso} then $f$ would be a degree five polynomial where all coefficients of the degree five terms are non-zero. However there is only one degree five term in $F$ so it is not possible. In case (ii) if $f$ has degree at least two then the polynomial had degree six at least since $g,h$ and $u$ have degree at least three. If $f$ is linear then there are no mixed terms in $F$ which is not the case here. The only remaining case is (iii), but note that $F(x,x,x)$ has four distinct roots, $\pm 1$ and $\pm \imath,$ while any $F(x,x,x)$ in case (iii) has at most three distinct roots.

To see that $F$ is a non-expander over the reals set $$A=\{\sqrt{2^i+1}\}_{i=1}^n, \; B=\{\sqrt{2^i+1}\}_{i=1}^n, \;  C=\{{2^i-1}\}_{i=1}^n.$$
With this choice of sets $A,B$ and $C,$ one can check by substitution that  $$|F(A,B,C)|\leq 3n.$$

\begin{example}
The polynomial $$F(x,y)=x^6+2x^4+2x^3y^3+4x^3+x^2+2xy^3+4x+y^6+4y^3+4$$
is an expander over the rationals.
\end{example}

\noindent
{\em Proof:} 
Before we analyze this example let us state Theorem \ref{vegso} for two variables. The proof is identical to the three-variable case.
\begin{theorem}\label{vegso2}
If $F(x,y)$ is not expander over the rationals then 
\begin{enumerate}[i]
\item , $F(x,y)=f(x+ay)$ or 
\item , $F(x,y)=f(g(x)+h(y))$ or 
\item , $F(x,y)=f((x+a)^{\alpha}(y+b)^{\beta})$
\end{enumerate}
for some (single variable) polynomial $f,$ polynomials $g, h$ with no constant terms which all have degree at least three and have at least two distinct roots, positive integers $\alpha,\beta$ and rational numbers $a,b.$ 
\end{theorem}
If $F$ had the form in (i) then all degree six terms were be there, however for example $x^5y$ has zero coefficient in $F.$ Let us suppose that $F$ has the form in (ii). $F$ has mixed terms, so $f$ has degree at least two and then it should be exactly two since $g,h$ have degree at least three. The lowest power of $y$ is three and $h$ has degree three too, so $h(y)=cy^3$ with some $c>0$ which is not possible because $h$ should have at least two distinct roots. The remaining case (iii) is not possible either since $F(x,x)$ has three distinct roots.

\section{Remarks}
Many questions remain open. 
\begin{itemize}
\item The most important problem for us would be to show that for any polynomial $f(x)$ over the rationals, $f(\mathbb{Q})$ contains no arbitrary long arithmetic progressions. (Or prove it for the integers at least) 
\item In our results above we proved that some rational polynomials -- which are non-expanders over the reals -- are expanders over the rationals. We believe that there are no integer polynomials which are expanders over the integers and non-expanders over the rationals, however we were unable to show it.
\item Our definition did not deal with the rate of expansion. It might be that the expanders over the rationals (or over the reals) are always good expanders; Is it true that for every expander $F(x,y)$ there are constants $c,\delta\geq 0$ so that $|F(A,A)|\geq c|A|^{1+\delta}$ for any $A\subset \mathbb{Q}$.  
\end{itemize}


\begin{thebibliography}{99}

\bibitem{AV}
D. Abramovich and J. F. Voloch, Lang's conjectures, fibered powers, and uniformity, New York J. of Math. II, pp 20--34,
1996. http://nyjm.albany.edu:8000/j/v2/Abramovich-Voloch.html


\bibitem{BGP} E. Bombieri, A. Granville, and J. Pintz,
Squares in arithmetic progressions,
Duke Mathematical Journal,  vol. 66 (1992)  165--204.

\bibitem{Bo}
J. Bourgain, 
More on the sum-product phenomenon in prime fields and its applications,
Int. J. Number Theory, 1(1), 1--32 (2005)


\bibitem{CHM} L. Caporaso, J. Harris, and B. Mazur
Uniformity of rational points, 
Journal of the 
AMS
Volume 10, Number 1, January 1997,  1--35.


\bibitem{Ch} M. C. Chang,
On Problems of Erd\H{o}s and Rudin. 
J. of Functional Analysis (2004) 43. 444--460. 


\bibitem{ER} Gy. Elekes and L. R\'onyai, 
A combinatorial problem on polynomials and rational functions, 
J. Combin. Theory, (A), 89 (2000), 1--20.

\bibitem{ESz} Gy. Elekes and E. Szab\'o,
How to find groups? (and how to use them in Erd\H{o}s geometry?),
Combinatorica, Online first (2012)
http://dx.doi.org/10.1007/s00493-012-2505-6

\bibitem{EL} Gy. Elekes
Sums versus products in Number Theory, Algebra and Erdšs Geometry --- a survey,
in: Paul Erd\H{o}s and his Mathematics II, Bolyai Math. Soc. Stud 11 Budapest, 2002, 241--290.

\bibitem{Fa} G. Faltings,  Endlichkeitss\"atze f\"ur abelsche Variet\"aten \"uber Zahlk\"orpern". Inventiones Mathematicae 73, 1983, (3): 349Ð366.

\bibitem{Fr}  G. A. Freiman, 
(1966) (in Russian). Foundations of a Structural Theory of Set Addition. 
Kazan: Kazan Gos. Ped. Inst.. pp. 140.


\bibitem{Pl} H. Pl\"unnecke, 
Eine zahlentheoretische anwendung der graphtheorie. J. Reine
Angew. Math. 1970, 243:171--183.

\bibitem{Ru} I. Z. Ruzsa,  
Generalized arithmetical progressions and sumsets. 
Acta Mathematica Hungarica (1994) 65 (4): 379--388


\bibitem{PR} I. Z. Ruzsa. An application of graph theory to additive number theory. Scientia,
Ser. A, 1989, 3: 97--109.

\bibitem{Sh} I. Shkredov, 
On monochromatic solutions of some nonlinear equations in $\mathbb{Z}_p$,
Mathematical Notes
2010, Volume 88, Issue 3-4,  603--611.



\bibitem{Si} J. H. Silverman,
A uniform bound for rational points on twists of a given curve,
J. London Math. Soc. (1993) vol. 47, 385--394.


\bibitem{SSZ}R. Schwartz, J. Solymosi, F. de Zeeuw
Extensions of a result of Elekes and R\'onyai,
arXiv:1206.2717 [math.CO]


\bibitem{So} J. Solymosi,
Arithmetic Progressions in Sets with Small Sumsets,
Combinatorics, Probability and Computing,
Volume 15 Issue 4, 2006, 597--603 
	arXiv:math/0503649 [math.NT]


\bibitem{Ta} T. Tao, 
Expanding polynomials over finite fields of large characteristic, and a regularity lemma for definable sets,
Cont. Disc. Math.(?)
arXiv:1211.2894 [math.CO]


\bibitem{TV} T. Tao and V. H. Vu, Additive Combinatorics,
Cambridge.


\bibitem{Xa} X. Xarles,
Squares in arithmetic progression over number fields.
Journal of Number Theory
Volume 132, Issue 3, 2012, Pages 379--389.
arXiv:0909.1642 [math.AG]


\end{thebibliography}
\end{document}